\documentclass[11pt]{amsart}

\usepackage{amssymb}

\usepackage{tikz}
\usetikzlibrary{patterns}
\usepackage[ruled,vlined,linesnumbered,norelsize]{algorithm2e}
\SetEndCharOfAlgoLine{}
\SetKwProg{Function}{Function}{}{}
\SetKwInOut{Input}{Input}
\SetKwInOut{Output}{Output}

\usepackage{hyperref} 
\usepackage{url}

\makeatletter
\newcommand{\addresseshere}{%
  \enddoc@text\let\enddoc@text\relax
}
\makeatother

\newcommand{\Z}{\mathbb{Z}}
\newcommand{\Q}{\mathbb{Q}}
\newcommand{\K}{\mathbf{K}}
\newcommand{\SL}{\mathrm{SL}}
\newcommand{\GFq}[1]{\mathbb{F}_{#1}}
\newcommand{\EE}{\mathcal{E}}
\newcommand{\Es}{\mathcal{E}^*}
\newcommand{\As}{A^*}
\newcommand{\Bs}{B^*}
\newcommand{\ComputeUMod}{{\sc ComputeUMod}}

\newcommand{\powersums}{\mathcal{P}}
\newcommand{\sfM}{{\sf M}}

\newtheorem{lemma}{Lemma}[section]
\newtheorem{theorem}[lemma]{Theorem}
\newtheorem{proposition}[lemma]{Proposition}

\begin{document}

\title{Computing the Charlap-Coley-Robbins \\ modular polynomials}
\author{François Morain}
\address{
    LIX - Laboratoire d'informatique de l'École polytechnique 
    \emph{and}
    GRACE - Inria Saclay--Île-de-France
}
\email{morain@lix.polytechnique.fr} 

\date{\today}

\maketitle

\begin{abstract}
Let $\EE$ be an elliptic curve over a field $\K$ and $\ell$ a prime.
There exists an elliptic curve $\Es$ related to $\EE$ by an
isogeny (rational map that is also a group homomorphisms) of
degree $\ell$ if and only $\Phi_\ell(X, j(\EE)) = 0$, where
$\Phi_\ell(X, Y)$ is the traditional modular polynomial.
Moreover, the modular polynomial
gives the coefficients of $\Es$, together with parameters needed to
build the isogeny explicitly. Since the traditional modular
polynomial has large coefficients, many families with smaller
coefficients can be used instead, as described by Elkies, Atkin and others.
In this work, we concentrate on the computation of the family of
modular polynomials introduced
by Charlap, Coley and Robbins. It has the advantage of giving directly
the coefficients of $\Es$ as roots of these polynomials. We review and
adapt the known algorithms to perform the computations of modular
polynomials. After describing the use of series computations, we
investigate fast algorithms using floating point numbers based on fast
numerical evaluation of Eisenstein series. We also explain how to use
isogeny volcanoes as an alternative.
\end{abstract}

%%%%% S
\section{Introduction}

Computing isogenies is the central ingredient of the
Schoof-Elkies-Atkin (SEA) algorithm
that computes the cardinality of elliptic curves over finite fields of
large characteristic~\cite{Schoof95,Atkin92b,Elkies98} and also
\cite{BlSeSm99}.
More recently, it has found its way in postquantum
cryptography~\cite{ChLaGo09,JaFe11,CaLaMaPaRe18,FeKoLePeWe20} among
others, as well as the cryptosystems~\cite{Couveignes06,RoSt06,FeKiSm18}.

Let $\K$ be a field of characteristic different from $2$ and $3$.
An isogeny between two elliptic curves $\EE/\K: Y^2 = X^3 + A X + B$ and
$\Es/\K: Y^2 = X^3 + \As X + \Bs$ is a group morphism that is a rational
map of degree $\ell$ (the cardinality of its kernel).
There are two ways to handle these isogenies. When the degree $\ell$
is small, formulas for $\As$, $\Bs$ and the kernel polynomial can be
precomputed. For large $\ell$,
one of the key ingredients is modular polynomials, the second one
finding rational expressions for $\As$ and $\Bs$ from $A$, $B$ and the
modular polynomial. We concentrate here on the former problem. The
second is treated in \cite{Morain23b}. Note there is a purely
algebraic approach using triangular sets~\cite[\S 7]{PoSc13}.

There are many families of modular polynomials that can be used, with
different properties. Very generally, a modular polynomial is some
bivariate polynomial $\Phi(X, J)$ where $J$ corresponds to the
$j$-invariant of the elliptic curve $\EE$, and $X$ stands for some
modular function on  $\Gamma_0(\ell)$.
The prototype is $j(\Es)$ (see below for
more precise statements) that
yields traditional modular polynomials. Alternative choices for $X$
exist. They all lead to polynomials of (conjectured) height $O((\ell+1)
\log\ell)$ but with small constants.

The Charlap-Coley-Robbins (CCR) modular polynomials~\cite{ChCoRo91}
offer an alternative as a triplet of polynomials $(U_\ell, V_\ell,
W_\ell)$ having the property (among others) that $\As$ (resp. $\Bs$)
is a root of $V_\ell$ (resp. $W_\ell$). They can be used in
conjunction with the preceding modular polynomials (in the case where
some partial derivative vanishes).

The aim of this work is describe the properties of CCR polynomials
(Section~\ref{sct:CCR}),
explain how to compute them (Section~\ref{sct:computations}) and use
them for isogeny computations. We
adapt methods from the classical approach: series
expansions over $\Z$ or $\GFq{p}$ for small primes and Chinese
remaindering theorem, evaluation/interpolation with floating point
numbers, isogeny volcanoes. For doing this, we need to evaluate
Eisenstein series (introduced in Section~\ref{sct:formulas}) at high
precisions and we describe algorithms to do so
in Section~\ref{sct:fast}, including an approach to the simultaneous
evaluation of several series. We give numerical examples and height
comparisons between modular polynomials.
We finish with the computation of
algebraic expressions for $\As$ and $\Bs$ as rational fractions
(see \cite{NoYaYo20}). 

%%%%% S
\section{Classical functions}
\label{sct:formulas}
\newcommand{\qt}{q_1}

Put $\qt = \exp(i \pi \tau)$ and $q = \qt^2$. Depending on authors,
formulas are expressed in either parameters, which sometimes is
clumsy. We write indifferently $f(\tau)$ or $f(q)$ some series.

%%%%%%%%%% SS
\subsection{Jacobi $\theta$ functions}

The classical $\theta$ functions are:
$$\theta_2(\qt) = \sum_{n\in\Z} \qt^{(n+1/2)^2},
\;\theta_3(\qt) = \sum_{n\in\Z} \qt^{n^2},
\;\theta_4(\qt) = \sum_{n\in\Z} (-1)^n \qt^{n^2}.$$
Among many properties, one has
$$\theta_{3}^4(\qt) = \theta_{4}^4(\qt) + \theta_{2}^4(\qt).$$
The latter formula enables to concentrate on the evaluation of
%% 10 = 2; 00 = 3; 01 = 4
$\theta_{3, 4}(\qt)$ as is done in \cite{Dupont11}.

Note also the following \cite[Prop. 4]{Dupont11}
\begin{proposition}\label{thinf}
$$\lim_{\mathrm{Im}(\tau) \rightarrow +\infty} \theta_{3}(\tau) = 1,
\; \lim_{\mathrm{Im}(\tau) \rightarrow +\infty} \theta_{4}(\tau) = 1,
\; \lim_{\mathrm{Im}(\tau) \rightarrow +\infty} \theta_{2}(\tau) = 0.$$
\end{proposition}

%%%%%%%%%% SS
\subsection{Eisenstein series}

%%%%%%%%%%%%%%% SSS
\subsubsection{Definitions, first properties}

The classical Eisenstein series\footnote{Ramanujan used $L = P = E_2$,
$M = Q = E_4$, $N = R = E_6$.} we consider are
$$E_2(q) = 1 -24 \sum_{n=1}^\infty \sigma_1(n) q^n,$$
$$E_4(q) = 1 + 240 \sum_{n=1}^\infty \sigma_3(n) q^n,$$
$$E_6(q) = 1 - 504 \sum_{n=1}^\infty \sigma_5(n)q^n,$$
where $\sigma_r(n)$ denotes the sum of the $r$-th powers of the
divisors of $n$. Other series $E_{2k}$ can be defined for $k > 3$.
Theory tells us that $E_{2k}$ is a modular form of weight $2k$ for $k
> 1$. As a result, a series for $k > 3$
can be expressed as polynomials in $E_4$ and $E_6$.
Also of interest is the discriminant $\Delta$:
$$\Delta(q) = (E_4(q)^3-E_6(q)^2)/1728 = \eta(q)^{24}$$
with $\eta(q) = q \prod_{n=1}^\infty (1-q^n)$ is the Dedekind
function. Finally, the modular invariant is
$$j(q) = \frac{E_4(q)^3}{\Delta(q)} = \frac{1}{q} + 744 + \cdots.$$

The quantities (see \cite[\S 13.20]{Erdelyi53})
$$a = \theta_{2}(\qt), b = \theta_{3}(\qt), c = \theta_{4}(\qt)$$
satisfy the following identities (among others)
\begin{equation}\label{E46D}
E_4 = (a^8+b^8+c^8)/2, \; E_6 = (a+b) (b+c) (c-a)/2, \; \Delta = (a b
c/2)^8.
\end{equation}
From which we deduce
$$\lim_{\mathrm{Im}(\tau) \rightarrow +\infty} E_4(\tau) = 1,
\; \lim_{\mathrm{Im}(\tau) \rightarrow +\infty} E_6(\tau) = 1.$$

%%%%%%%%%%%%%%% SSS
\subsubsection{The special case of $E_2$}

The series $E_2$ is not a modular form since (see \cite{Rankin77}):
\begin{theorem}
For all matrices
$\left(\begin{array}{cc} a & b \\ c & d \\\end{array}\right)$
of $\SL_2(\Z)$, one has
$$E_2((a\tau+b)/(c\tau+d)) = (c\tau+d)^2 E_2(\tau) + \frac{6 c}{\pi i} (c \tau
+ d).$$
\end{theorem}
We can build a modular form easily as follows, using
{Ramanujan's multipliers}. Let $n$ be an integer and
Let $\mathcal{F}_n$ denote the {\em multiplier} $E_2(\tau)-n
E_2(n\tau)$.
\begin{proposition}
The function $\mathcal{F}_n$ is a modular form of weight 2 and
trivial multiplier system for $\Gamma_0(n)$.
\end{proposition}

\noindent
{\em Proof:} Write, for $ad-bc=1$ the value
$$E_2\left(n \frac{a\tau+b}{c\tau+d}\right)
= E_2\left(\frac{a (n\tau)+n b}{(c/n)(n\tau)+d}\right)$$
$$= ((c/n)n\tau+d)^2 E_2(n\tau) - \frac{6 c}{n \pi i}
((c/n) n\tau + d)$$
which leads to
$$n E_2\left(n \frac{a\tau+b}{c\tau+d}\right)
= n (c\tau+d)^2 E_2(n\tau) - \frac{6 c}{\pi i}
(c\tau + d).$$
Subtracting $E_2((a\tau+b)/(c\tau+d))$, we see that
$$\mathcal{F}_n((a\tau+b)/(c\tau+d)) = (c\tau+d)^2
\mathcal{F}_n(\tau). \Box$$

Some identities are known for small values of $n$, for instance
\cite{KaKo03} for $n \in \{2, 4\}$; \cite{Berndt89} for $n=3$ and $11$;
\cite[Thm 6.2]{BoBo91}, 
\cite[Thm 3.7]{BeChSoSo00} for $n \in \{5, 7\}$. A very nice
relation is~\cite[Thm 6.3]{BeChSoSo00}
$$\mathcal{F}_7(q) = 6 \left(\sum_{m, n = -\infty}^\infty
q^{m^2+mn+2n^2}\right)^2.$$

%%%%% S
\section{Fast evaluations}
\label{sct:fast}

%%%%%%%%%% SS
\subsection{Fast evaluation of $E_{2k}(q)$ for $k \geq 1$}

The $\theta$ functions that can be evaluated at precision $N$ in
time $O({\sf M}(N) \sqrt{N})$ with $\qt$-expansions (see
\cite{EnHaJo18}) or faster in $O({\sf M}(N)
\log N)$ using \cite{Dupont11} and also \cite{Labrande18}.
It follows that the quantities $E_{2k}$ (for $k \geq 2$)
can be evaluated at precision $N$ in $O({\sf M}(N)\log N)$
operations. As a consequence $j(q)$ can also be evaluated with the
same complexity.

\newcommand{\hypergeom}[4]{{}_2F_1\left({#1}, {#2}; {#3}; {#4}\right)}
Evaluating $E_2$ is less obvious. However, hidden in the proof of
\cite[Thm 4]{KaZa98} (thanks to \cite{KaKo03} for highlighting this),
we find
$$\frac{E_2 E_4}{E_6} = 
\frac{\hypergeom{\frac{13}{12}}{\frac{5}{12}}{1}{\frac{1728}{j}}}
     {\hypergeom{\frac{1}{12}}{\frac{5}{12}}{1}{\frac{1728}{j}}}
   = 1 + \frac{720}{j} + \cdots$$
where the Gauss hypergeometric function is defined by
$$\hypergeom{a}{b}{c}{x} = \sum_{k=0}^\infty \frac{(a)_k (b)_k}{(c)_k \,
k!} x^k, |x| < 1$$
where $(a)_k = a (a+1) \cdots (a+k-1)$. By
   \cite{Hoeven99,Hoeven01,MeSa10} and also \cite{BrZi10},
   this function can be computed at precision $N$ in $O({\sf
   M}(N)(\log N)^2)$ operations. See also \cite{Johansson19} for
   realistic computations.

Other links with hypergeometric
functions could be investigated (A~.Bostan, personal communication).

Also, note that actually we need to
evaluate $\mathcal{F}_\ell$ for prime $\ell$. Sometimes, we may use
special formulas as indicated above. 

%%%%%%%%%% SS
\subsection{A multi-value approach}

In practice, a simpler approach yields the values $E_{2k}$ of many
$k$'s with $k \geq 1$ in time $O({\sf M}(N) \sqrt{N})$ based on
\cite{EnHaJo18}. The cost reduces to that of one series evaluation.

From \cite{BeYe02}, we take
$$(q; q)_\infty = \exp(-2i\pi\tau/24) \eta(q),
$$
and for $k \geq 0$:
$$T_{2k}(q) = 1 + \sum_{n=1}^\infty (-1)^n \left\{(6n-1)^{2k} q^{n
(3n-1)/2}\right.$$
$$\left.+ (6n+1)^{2k} q^{n (3n+1)/2}\right\}.$$
Note that $(q; q)_\infty = T_0(q)$.

\begin{theorem}[Section 6]\label{ET}
$$\frac{T_2(q)}{T_0(q)} = E_2,
\;\frac{T_4(q)}{T_0(q)} = 3 E_2^2-2 E_4,
\;\frac{T_6(q)}{T_0(q)} = 15 E_2^3 - 30 E_2 E_4 + 16 E_6.$$
\end{theorem}
General formulas for $T_{2k}(q)/T_0(q)$ are given in the reference.

If we need to compute $E_2$, $E_4$ and $E_6$, we see that it is enough to
evaluate the series $T_{2k}$ for $k \in \{0, 1, 2, 3\}$ followed by a
handful of multiplications and divisions as given in the preceding
Theorem. Moreover, we can evaluate these series by
sharing the common powers of $q$. These powers are evaluated at a
reduced cost using \cite[Algorithm2]{EnHaJo18}. We give the modified
procedure as algorithm~\ref{algo1}.
In Step 3.2.3, we have added the contribution $(6n \pm 1)^{2i}$ to each
\verb+T[i]+. We assume that the cost of multiplying by these small
quantities is negligible. Were it not the case, we could use
incremental computations of the polynomials $(6n \pm 1)^{2i}$. The
cost of this algorithm reduces to that of one of the series, gaining a
factor $kmax$.

\LinesNotNumbered
\begin{algorithm}[hbt]
\caption{Combined evaluation of $T_{2k}(q)$. \label{algo1}}
\SetKwProg{Fn}{Function}{}{}
\Fn{EvaluateManyT(q, N, kmax)}{
\Input{$q$, $N$, $kmax$}
\Output{$(T_{2k}(q))$ for $0 \leq k \leq kmax$}

1. \For{$k:=0$ \KwTo $kmax$}{
    $T[k] \leftarrow 0$\;
}

2. $s \leftarrow 1$; $A \leftarrow \{1\}$; $Q[1] \leftarrow \{q\}$; $c
\leftarrow 0$\;

3. \For{$n:=1$ while $n (3n+1)/2 \leq N$}{
$s \leftarrow -s$\; \tcp*[h]{$s = (-1)^n$}

3.1 $c \leftarrow c + 2 n-1$\;

3.2 \For{$r:=1$ \KwTo $2$}{
 3.2.1 \If{$r = 2$}{
    $c \leftarrow c + n$\; \tcp*[h]{$c = n (3n+1)/2$}
 }
 3.2.2 $q' \leftarrow $ FindPowerInTable($A$, $Q$, $c$)\;
 3.2.3 $C \leftarrow (6n+(-1)^r)^2$\;
 3.2.4 \For{$k:=0$ \KwTo kmax}{
     $T[k] \leftarrow T[k] + s q'$\;
     \If{$k < kmax$}{
         $q' \leftarrow C q'$\;
     }
 }}}
4. \For{$k:=0$ \KwTo kmax}{
   $T[k] \leftarrow T[k] + 1$\;
   }
5. \Return{$T$}.
}
\end{algorithm}
Algorithm~\ref{algo1} uses the primitive in Algorithm~\ref{algo0}. The
reason of Step 4 is that $T[k]$ will be close to $1$ when $q$ is
small, so that we may not want to add 1 right at the beginning and
perhaps not in this function.

\LinesNotNumbered
\begin{algorithm}[hbt]
\caption{Finding $c$ as a combination of known values. \label{algo0}}
\SetKwProg{Fn}{Function}{}{}
\Fn{FindPowerInTable($A$, $Q$, $c$)}{
\Input{$A = \{a_1, \ldots, a_z\}$, $Q$ such that for all $i$, $Q[a_i]
= q^{a_i}$, $c$}
\Output{$q^c$; $A$ and $Q$ are updated}

 \If{$c = 1$}{
    $q' \leftarrow Q[1]$\;
 }
 \ElseIf{$c = 2a$ with $a \in A$}{
    $q' \leftarrow Q[a]^2$\;
 }
 \ElseIf{$c = a+b$ with $a, b \in A$}{
    $q' \leftarrow Q[a] \cdot Q[b]$;
 }
 \ElseIf{$c = 2 a+b$ with $a, b \in A$}{
    $q' \leftarrow Q[a]^2 \cdot Q[b]$\;
 }
   $A \leftarrow A \cup \{c\}$\;
   $Q[c] \leftarrow q'$\;
   \Return{$q'$}.
}
\end{algorithm}

%%%%%%%%%% SS
\subsection{The case of imaginary arguments}

In practice, it is easier to consider $\tau = \rho i$ for real $\rho \geq
1$. In that case, $1 > q_0 = \exp(-2\pi) \geq q = \exp(-2\pi \rho) >
0$. The functions $E_2$ and $E_6$ are increasing from $E_{2k}(q_0)$ to
$1$ (note that $E_6(q_0) = 0$ and $E_2(q_0) = 3/\pi$ from~\cite{ElSe10});
$E_4$ is decreasing from $E_4(q_0)$ to
$1$. This is important to note for the computations not to
explode. Remember also that $j(i) = 1728$.

We turn to the precision needed for evaluating the functions $T_{2k}$.
Let $N$ denote an integer and $T_{2k, N}$ the truncated sum up to $n =
N-1$. Since the series is alternating, we can bound the error using

$|T_{2k}(q)-T_{2k, N}(q)|$

$\leq \{(6N-1)^{2k} \, q^{N(3 N-1)/2} +
(6N+1)^{2k} \, q^{N(3 N+1)/2}\}$

$\leq ((6N-1)^{2k}+(6N+1)^{2k}) \, q^{N(3 N-1)/2}.$

\noindent
Since $0 < q < 1$, this gives us a quadratic convergent series.

%%%%% S
\section{The polynomials of Charlap-Coley-Robbins}
\label{sct:CCR}

%%%%%%%%%% SS
\subsection{Division polynomials}

For $\EE: y^2=x^3+A x+B$, multiplication of a point $(X, Y)$ by positive
$n$ on $\EE$ is given by
$$[n] (X, Y) = \left(\frac{\phi_n(X, Y)}{\psi_n(X, Y)^2}, \frac{\omega_n(X,
Y)}{\psi_n(X, Y)^3}\right)$$
where the polynomials satisfy
$$\phi_n = X \psi_n^2 - \psi_{n+1} \psi_{n-1},
\; 4 Y \omega_n = \psi_{n+2} \psi_{n-1}^2 -
\psi_{n-2}\psi_{n+1}^2$$
and $\phi_n, \psi_{2n+1}, \psi_{2n}/(2Y), \omega_{2n+1}/Y,
\omega_{2n}$ belong to $\Z[A, B, X]$. It is customary to simplify this
using
$$f_n(X) = \left\{\begin{array}{ll}
\psi_n(X, Y) & \mathrm{for}\; n \;\mathrm{odd}\\
\psi_n(X, Y)/(2 Y) & \mathrm{for}\; n \;\mathrm{even}
\end{array}\right.$$
with first values
$$f_{-1} = -1, \; f_0 = 0, \; f_1 = 1, \; f_2 = 1,$$
$$f_3(X, Y) = 3 X^4 + 6 A X^2 + 12 B X - A^2$$
$$f_4(X, Y) = X^6 + 5 A X^4 + 20 B X^3 - 5 A^2 X^2
- 4 A B X -8 B^2 - A^3.$$
The degree of $f_n$ is $(n^2-1)/2$ for odd $n$ and $(n^2-4)/2$ for
even $n$. If $X$ has weight $1$, $A$ weight 2 and $B$ weight 3, all
monomials in $f_n$ have the same weighted degree equal to the degree
of $f_n$.

%%%%%%%%%% SS
\subsection{The work of Elkies}

An isogeny is associated with its kernel, or its polynomial
description (called {\em kernel polynomial}). Given some finite
subgroup $F$ of $\EE$, one can build an
isogenous curve $\Es$ and the corresponding isogeny, using V\'elu's formulas.
In the context of point counting, we discover
a curve $\Es$ that is $\ell$-isogenous to $\EE$ via its $j$-invariant
as a root of the traditional
modular polynomial, and we need to find the coefficients of $\Es$,
together with the isogeny. The idea of Elkies is to consider the same
problems on the Tate curves associated to the elliptic curves $\EE$
and $\Es$.

To be brief, $\EE$ has an equation in some parameter $q$, and the
isogenous $\Es$ is associated to parameter $q^\ell$, where $\ell$ is
the degree of the isogeny, which in our case is associated with a
finite subgroup $F$ of cardinality $\ell$. To be more precise, we
consider $\EE$ has having equation $Y^2 = X^3 + A X + B$ with
$$A = -3 E_4(q), B = -2 E_6(q).$$
With a compatible scaling, we get the equation for $\Es: Y^2 = X^3 +
\As X + \Bs$ with
$$\As = -3 \ell^4 E_4(q^\ell), \; \Bs = -2 \ell^6 E_6(q^\ell).$$
More importantly, writing $\sigma_r$ for the power sums of the roots of
the kernel polynomial, we have
$$\sigma_1 = \frac{\ell}{2} (\ell E_2(q^\ell)-E_2(q)) =
-\frac{\ell}{2} \mathcal{F}_\ell(q).$$
Beyond this, Elkies proved~\cite[formulas (66) to (69)]{Elkies98}
\begin{proposition}\label{Elkies:formulas}
$$A-\As = 5 (6 \sigma_2 + 2 A \sigma_0),$$
$$B-\Bs = 7 (10 \sigma_3 + 6 A \sigma_1+4 B \sigma_0),$$
together with an induction relation satisfied by other $\sigma_k$ for
$k > 3$.
\end{proposition}
This can rephrased as $(\sigma_1, \As, \Bs)$ is enough to describe an
isogeny. Also $\As$ and $\Bs$ belong to $\Q[\sigma_1, A, B]$
since $\sigma_2$ and $\sigma_3$ do. The minimal polynomial of
$\sigma_1$ is a modular polynomial, and we can express $\As$ and $\Bs$
as elements in the field $\Q[\sigma_1, A, B]$, which we use
below. Rephrased another times, $E_4(q^\ell)$ and $E_6(q^\ell)$ are
modular forms we need to express as known modular
forms. See~\cite{Elkies98} for more details on this subject.

Given these quantities, there are several algorithms to get the
isogeny. We refer to \cite{BoMoSaSc08} for this.

%%%%%%%%%% SS
\subsection{The CCR polynomials}

%%%%%%%%%%%%%%% SSS
\subsubsection{Reinterpreting Elkies's results}

One way of looking at the work of Elkies (taken from
\cite{ChCoRo91}) is to realize that we try to decompose the polynomial
$f_\ell$ (say $\ell$ is odd) of degree $(\ell^2-1)/2$ over a subfield
of degree $\ell+1$.

\begin{center}
\begin{tikzpicture}[scale=0.5]
\node (a) at (2, 2) {$\Q(A, B)[X]/(f_\ell(X, A, B))$};
\node (b) at (2, 0) {$\Q(A, B)[X]/(U_\ell(X, A, B))$};
\node (c) at (2, -2) {$\Q(A, B)$};
\draw[-] (a) -- (b) node[midway,right] {$(\ell-1)/2$};
\draw[-] (b) -- (c) node[midway,right] {$\ell+1$};
\end{tikzpicture}
\end{center}

A classical way for doing this is to use the
trace $t_1$ of an element in $\Q(A, B)[X]/(f_\ell(X, A, B))$. Let $x_1$
stand for the (formal) abscissa of an $\ell$-division point $P = (x_1,
y_1) \neq O_E$. Other points are $P_j = [j] P = (x_j, y_j)$ and $x_j$
can be expressed using division polynomials. For $0 \leq k \leq
\ell+1$, we define
\begin{equation}\label{defpk}
t_k = \sum_{j=1}^d x_j^k = \sum_{j=1}^d \left(x_1 - \frac{\psi_{j-1}(x_1)
\psi_{j+1}(x_1)}{\psi_j(x_1)^2} \right)^k
\end{equation}
so that $t_1 = x_1 + \cdots + x_d$ and $t_0 = d = (\ell-1)/2$. The
minimal polynomial $U_\ell(X) = X^{\ell+1} + u_1 X^{\ell} +
\cdots + u_0$ of $t_1$ defines the lower     
extension. Given the power sums $t_k$'s, Newton's identities enable us
to reconstruct the minimal polynomial $\prod_{i=1}^d (X-x_i) = X^d -
t_1 X^{d-1} + \cdots$ over the intermediate extension.
Putting everything together, $t_1$ coincides with $\sigma_1$ included
above.

By direct application of (\ref{defpk}), we can compute $U_3(X, Y, Z)
= X^4+2 Y X^2+4 Z X-Y^2/3$ (this is $\psi_3$ in disguised
form). Larger values of $\ell$ require more work.

%%%%%%%%%% SS
\subsubsection{Theory}

We start from an elliptic curve $\EE: Y^2 = X^3+A X+B$ and we fix some
odd prime $\ell$, putting $d = (\ell-1)/2$. Our aim is to find the
equation of an $\ell$-isogenous curve $\Es: Y^2 = X^3 + \As X + \Bs$. 

\begin{theorem}
There exist three polynomials $U_\ell$, $V_\ell$,
$W_\ell$ in $\Z[X, Y, Z, 1/\ell]$ of degree $\ell+1$ in $X$
such that $U_\ell(\sigma_1, A, B)=0$, respectively $V_\ell(\As, A, B) = 0$,
$W_\ell(\Bs, A, B) = 0$.
\end{theorem}

Let us turn our attention to the properties of these polynomials.
\begin{theorem}
When $\ell > 3$, the polynomials $U_\ell$, $V_\ell$, $W_\ell$ live in
$\Z[X, Y, Z]$.
\end{theorem}

\begin{proposition}
Assigning respective weights 1, 2, 3 to $X$, $Y$, $Z$,
the polynomials $U_\ell$, $V_\ell$ and $W_\ell$ are homogeneous with weight
$\ell+1$.
\end{proposition}

\begin{proposition}
The roots of $U_\ell(X, A(q), B(q))$ are $-\ell
\mathcal{F}_\ell(q)/2$ and $\mathcal{F}_\ell(w \zeta_\ell^k)/2$ for $0
\leq k < \ell$, where $w^\ell = q$ and $\zeta_\ell$ is a root of unity.
\end{proposition}

We follow \cite[\S 8]{ChCoRo91}.
\begin{proposition}\label{sizeUVW}
The height of $U_\ell$ (resp. $V_\ell$, $W_\ell$) is approximately $2
k (\ell+1) \log\ell$ for $k=1, 2, 3$ corresponding to $U$, $V$, $W$
respectively.
\end{proposition}

\medskip
\noindent
{\em Proof:} let $\lambda_k(u) =
\sum_{n=1}^\infty \delta_{2k-1}(n) u^n$ which the important term in
the sums. Using 
$$\delta_{2k-1}(n) \leq n^{2k-1} \sum_{d \mid n} \frac{1}{d^{2k-1}}
\leq \zeta(2k-1) n^{2k-1},$$
we can see that $\lambda_k(u)$ is dominated by $\zeta(2k-1)
\sum_{n=1}^\infty n^{2k-1} u^n$ dominated by
$$(2k-1)! \, \zeta(2k-1) \sum_{n=1}^\infty \binom{n+2k-2}{2k-1} u^n = C_k u
(1-u)^{-2k}.$$
So $\lambda_k^m$ is dominated by $C_k^m (1-u)^{-2km}$, and the
highest term is approximately $u^{m^2/2}$. An approximation to the
largest coefficient of $\lambda^m$ is 
$$C_k^m \binom{m^2/2+2k m -1}{2km-1} \approx C_k^m
\frac{(m^2/2+2km)}{(2km/e)}^{2km}$$
$$ = C_k^m \left(\frac{e}{4k}\right)^{2k
m} (m + 4k)^{2k m}.$$
Taking $m = \ell+1$ leads to the result. $\Box$

The traditional modular polynomial $\Phi_\ell^t$ has height $6
(\ell+1)\log\ell$ approximately~\cite{Cohen84,Sutherland13}. We see
that $U_\ell$ and $V_\ell$ are smaller, but that $W_\ell$ is as big as
$\Phi_\ell^t$. See the appendix for some tables.

%%%%%%%%%% SS
\subsubsection{Computing isogenous curves over finite fields}
\label{sct:summary}

When using $U_\ell$, $V_\ell$, $W_\ell$, we need to find the roots of
three polynomials of degree $\ell+1$ instead of $1$ in the traditional
case. In general, if
$U_\ell$ has rational roots (it should be 1, $2$ or $\ell+1$), then
this is the case for each of $V_\ell$, $W_\ell$. For each triplet of
solutions $(\sigma_1, z_1, z_2)$ we need to test whether this leads to
an isogeny or not. See techniques for this task in \cite{BoMoSaSc08}.

%%%%% S
\section{Computing CCR polynomials}
\label{sct:computations}

Using the results of the preceding section, the polynomials can be
written as
$$U_\ell(X, Y, Z) = X^{\ell+1} + \sum_{i_1+2i_2+3i_3 = \ell+1}
u_{i_1, i_2, i_3} X^{i_1} Y^{i_2} Z^{i_3},$$
$$V_\ell(X, Y, Z) = X^{\ell+1} + \sum_{2 i_1+2i_2+3i_3 = \ell+1}
v_{i_1, i_2, i_3} X^{i_1} Y^{i_2} Z^{i_3},$$
$$W_\ell(X, Y, Z) = X^{\ell+1} + \sum_{3 i_1+2i_2+3i_3 = \ell+1}
w_{i_1, i_2, i_3} X^{i_1} Y^{i_2} Z^{i_3}.$$
All the methods to be described can be applied to $U_\ell$, $V_\ell$
and $W_\ell$. To simplify the presentation,
we assume from now on that $\ell > 3$ and concentrate on $U_\ell$,
indicating what has to be changed for $V_\ell$ (resp. $W_\ell$).
We rewrite
$$U_\ell(X, Y, Z) = X^{\ell+1} + \sum_{r=0}^{\ell}
X^r \sum_{2i_2+3i_3 = \ell+1-r} c_{r, i_2, i_3} Y^{i_2} Z^{i_3}.$$
From this, we can see that there are no possible terms for $r =
\ell$.

We first count the number of monomials.
We take the following from \cite[p. 110]{Comtet74}.
\begin{proposition}\label{propN3}
The number of solutions $\mathbf{N}_{1, 2, 3}(n)$ in positive integers
of $i_1 + 2 i_2 + 3 i_3 = n$ is the closest integer to $(n+3)^2/12$.
\end{proposition}
Using the same method
\begin{proposition}\label{propN2}
The number of solutions $\mathbf{N}_{2, 3}(n)$ in positive integers of $2
i_2 + 3 i_3 = n$ is
$$\mathbf{N}_{2, 3}(n)
= \frac{n+1}{6}+\frac{(-1)^n}{4}-\frac{1}{12}+
\left\{\begin{array}{cc}
\frac{2}{3} & \text{ if } n \equiv 0 \bmod 3,\\
\frac{1}{3} & \text{ if } n \equiv 2 \bmod 3.\\
\end{array}\right.$$
\end{proposition}
We are going to compute many products of the form $Y^{i_2} Z^{i_3}$
for $Y$ and $Z$ that are large floating point numbers or
series. We use Pippenger's algorithm~\cite[pp. 247--249]{Pippenger80}
(thanks to~\cite{Bernstein02b}) for that task.

The authors of \cite{ChCoRo91} give three methods to compute the
polynomial $U_\ell$. The first is based on equation (\ref{defpk})
and can be used for very small $\ell$'s.
Two more methods use manipulations of $q$-expansions of series over
$\Q$, or modulo small primes followed by recovery using the Chinese
remaindering theorem using the bounds in Proposition~\ref{sizeUVW}.

%%%%%%%%%%%%%%%	SSS
\subsection{Using $q$-expansions}

Note that
$$\sigma_1(q) = -\frac{\ell\mathcal{F}_\ell(q)}{2} = \frac{\ell
(\ell-1)}{2} + 12 \ell \sum_{n=1}^\infty \sigma_1'(n) q^n$$
where $\sigma_1'(n)$ is the sum of the divisors of $n$ prime to $\ell$.

The second method proceeds
by plugging the series $\sigma_1(q)$, $A(q)$ and $B(q)$ up to degree
$\mathbf{N}_{1, 2, 3}(\ell+1) = O(\ell^2)$ (Proposition \ref{propN3})
in $U_\ell(X, Y, Z) = 0$
and find the unknown coefficients of the polynomial using a linear
system over the rationals (in fact integers for $\ell > 3$). This can
be done for small $\ell$'s using any mathematical system. Note that
the number of algebraic operations (multiplications in $\Q$, or a
finite field) will be close to $O(\ell^{2\omega})$ where $\omega$ is
the constant for matrix multiplication with $2 \leq \omega \leq 3$, so
typically $O(\ell^6)$, which is large.

The third method exploits the fact that the power sums
$\sigma_r(q)$ for $1 \leq r \leq \ell+1$ are modular forms and can be
represented as polynomials in $A(q)$ and $B(q)$ (or $E_4(q)$, $E_6(q)$)
$$\sigma_r(q) = \sum_{2 i_2 + 3 i_3 = r} u_{r, i_2, i_3} A(q)^{i_2}
B(q)^{i_3}.$$
This leads to a linear system $\mathcal{S}_r$ (independent of $\ell$)
in the $u_{r, i_2, i_3}$'s. By Proposition~\ref{propN2}, the system
has size $\mathbf{N}_{2, 3}(r) \times \mathbf{N}_{2, 3}(r) \approx
(r/6)^2$ and the linear system can be solved with $O(r^{\omega})$
operations over $\Z$, for a total of $O(\sum_r r^{\omega}) =
O(\ell^{\omega+1}) = O(\ell^4)$ generally. Note that all these systems
may be solved in parallel.

The authors of \cite{ChCoRo91} suggest to use a more adapted
basis. When $r = 2m$, use $\{E_4^{m-3j} \Delta^k, 0 \leq k \leq
\lfloor m/3\rfloor\}$; when $r = 2m+3$ (remember that the coefficient
for $r = 1$ is 0), use $\{E_6 E_4^{m-3k} \Delta^k, 0 \leq k \leq
\lfloor m/3\rfloor\}$. Note that in all cases, the series for index $k$
start with $q^k$ and the bases present a triangular shape, in other
words the corresponding system $\mathcal{S}_r$ is
triangular. Moreover, since the leading coefficient is 1, all
solutions to the system will be integers.

Once solved for all $r$'s, we use
Newton's identities to recover the coefficients of $U_\ell$. One needs
to evaluate the series $\sigma_r(q)$
using intermediate expressions in $w = q^{1/\ell}$
having roots of unity $\zeta_\ell$ temporarily appearing and vanishing.
See \cite[\S 2.2]{Enge09b} for more details and complexity analysis.
In particular, if we denote by $\sfM_q(d)$ the number of arithmetic
operations in $\Z$ required to multiply two dense $q$-expansions with
$d$ terms, then the total complexity of the series computations is
$O(\ell \sfM_q(\ell d)$, which is $O(\ell \sfM_q(\ell^2))$ in our
case. If $H$ is a bound on the height of the polynomial, then the bit
complexity is $O(\ell^3 (\log \ell) \sfM(H))$. Assuming $H \in
O(\ell\log\ell)$ by Proposition~\ref{sizeUVW}, this is $O(\ell^4
\log^{3+\epsilon} \ell)$.

\medskip
\noindent
{\bf Example.} Consider the case $\ell = 5$. The systems
$\mathcal{S}_r$ to be solved come from the equations:
\begin{eqnarray*}
\sigma_2(q) &=& u_{2, 0} E_4 \\
\sigma_3(q) &=& u_{3, 0} E_6 \\
\sigma_4(q) &=& u_{4, 0} E_4^2 \\
\sigma_5(q) &=& u_{5, 0} E_6 E_4 \\
\sigma_6(q) &=& u_{6, 0} E_4^3 + u_{6, 1} \Delta
\end{eqnarray*}
We compute
$$\sigma_6(q) = 1000320+186071040 q+ \cdots$$
and we remember that $E_4(q) = 1+240 q + \cdots$, $\Delta(q) = q + \cdots$ so
that the system $\mathcal{S}_6$ is
$$\left\{\begin{array}{rcl}
1000320 &=& u_{6, 0} \\
186071040 &=& 720 \, u_{6, 0} + u_{6, 1}\\
\end{array}\right.$$
which is triangular indeed and therefore easy to solve. Its solutions
are integers.

We can also work modulo small primes and use the Chinese Remainder
Theorem to recover the polynomials. 

%%%%%%%%%% SS
\subsection{Floating point methods}

We adapt the methods proposed for ordinary modular equations to our
triplet of polynomials $(U_\ell, V_\ell, W_\ell)$. We note $H$ for the
logarithmic height of the polynomials, that we have estimated to $2k
(\ell+1)\log\ell$ in Proposition~\ref{sizeUVW}.

%%%%%%%%%%%%%%% SSS
\subsubsection{Solving a linear system}

We start from
$$U_\ell(\sigma_1, A, B) = 0 = \sigma_1(q)^{\ell+1} $$
$$+ \sum_{r=0}^{\ell}
\sigma_1(q)^r \sum_{2 i_2 + 3 i_3 = \ell+1-r} c_{r, i_2, i_3} A^{i_2}(q)
B^{i_3}(q)$$
and we compute floating point values to get a linear system in the
$c_{r, i_2, i_3}$ that should come out as integers for $\ell > 3$. We
evaluate $\sigma_1(q)$, $A(q) = -3 E_4(q)$ and $B(q) = -2 E_6(q)$ at
high precision for chosen values of $\tau$ in $q =
\exp(2i\pi\tau)$. This would involve $O(\ell^{2\omega})$ floating
point operations, and we can do better in the next section.

%%%%%%%%%%%%%%% SSS
\subsubsection{Using power sums}

The case of the traditional modular polynomial is treated in
\cite{Enge09b}. We can use the same approach for our polynomials.
First of all, we need to compute
$$-\frac{\ell}{2} \mathcal{F}_\ell(q), \{\mathcal{F}_\ell(w
\zeta_\ell^k)/2, 0 \leq k < \ell\}.$$
By definition
$$\mathcal{F}_\ell(q) = E_2(q)-\ell E_2(q^\ell)$$
and
$$\mathcal{F}_\ell(w \zeta_\ell^k) = E_2(w \zeta_\ell^k) - \ell
E_2(q),$$
and the last term is a constant w.r.t. $k$. We first evaluate
$E_2(q)$, $E_2(q^\ell)$ and then the other roots, by sharing the
computations: All terms we need are of the form $(w \zeta_\ell^k)^e = w^e
\zeta_\ell^{(k e) \bmod \ell}$. When $\ell \mid e$, the computation is
a little faster. We give the corresponding code as
Algorithm~\ref{algo2}. We also precompute $\xi_k = \zeta_\ell^k$. The
complete code is in Algorithm~\ref{algo3}.

\LinesNotNumbered
\begin{algorithm}[hbt]
\caption{Combined evaluation of $T_{2k}(w \zeta_\ell^j)$. \label{algo2}}
\SetKwProg{Fn}{Function}{}{}
\Fn{EvaluateConjugateValues($\ell$, $w$, $N$, $(\xi)$, $kmax$)}{
\Input{$\ell$, $w$, $(\xi)$, $N$}
\Output{$(T_{2k}(w \zeta_\ell^j))$ for $0 \leq k \leq kmax$, $0 \leq j
< \ell$}

1. \For{$k:=0$ \KwTo $kmax$}{
     \For{$j:=0$ \KwTo $\ell-1$}{
         $T[k, j] \leftarrow 0$\;
     }
}

2. $s \leftarrow 1$; $A \leftarrow \{1\}$; $W[1] \leftarrow \{w\}$; $c
\leftarrow 0$\;

3. \For{$n:=1$ while $n (3n+1)/2 \leq N$}{
$s \leftarrow -s$\; \tcp*[h]{$s = (-1)^n$}

3.1 $c \leftarrow c + 2 n-1$\;

3.2 \For{$r:=1$ \KwTo $2$}{
 3.2.1 \If{$r = 2$}{
    $c \leftarrow c + n$\; \tcp*[h]{$c = n (3n+1)/2$}
 }
 3.2.2 $w' \leftarrow s \cdot$FindPowerInTable($A$, $W$, $c$)\;
 3.2.3 $C \leftarrow (6 n+(-1)^r)^2$\;
 3.2.4 \For{$k:=0$ \KwTo kmax}{
     $T[k, 0] \leftarrow T[k, 0] + w'$\;
     \For{$j:=1$ \KwTo $\ell-1$}{
         $T[k, j] \leftarrow T[k, j] + \xi[(j c) \bmod \ell] \, w'$\;
     }
     \If{$k < kmax$}{
         $w' \leftarrow C \cdot w'$\;
     }
}}}
4. \For{$k:=0$ \KwTo $kmax$}{
     \For{$j:=0$ \KwTo $\ell-1$}{
         $T[k, j] \leftarrow T[k, j] + 1$\;
     }
   }
5. \Return{$T$}.
}
\end{algorithm}

\LinesNotNumbered
\begin{algorithm}[hbt]
\caption{Computing CCR polynomial $k$ using floating point
 numbers. \label{algo3}}
\SetKwProg{Fn}{Function}{}{}
\Fn{ComputeUVW($k$, $\ell$)}{
\Input{$k \in \{1, 2, 3\}$ corresponding to $U_\ell$, $V_\ell$ or
$W_\ell$ respectively, $\ell$ an odd prime}
\Output{the CCR polynomial of order $k$}

0.0 $H \leftarrow 2k (\ell+1) \log \ell$; all computations are carried
out at precision $H$\;

0.1 compute $\zeta_\ell \leftarrow \exp(2 i \pi/\ell)$\;

0.2 \For{$j:=0$ \KwTo $\ell-1$}{
    $\xi_j \leftarrow \zeta_\ell^j$\;
   }

0.3 Compute all systems $\mathcal{S}_r$ for $2 \leq r \leq \ell+1$\;

0.4 \For{$r:=2$ \KwTo $\ell+1$}{
    $\mathcal{L}_r \leftarrow \emptyset$\;
    }

0.5 $\rho \leftarrow 1$\;

1. \While{there is a system $\mathcal{L}_r$ that is not solved}{

1.0 $\rho \leftarrow \rho + 0.1$\;

1.1 $w \leftarrow \exp(-2 \pi \rho/\ell)$\;

1.2 $T \leftarrow $ EvaluateConjugateValues($\ell$, $w$, $N$, $(\xi)$,
$2k$)\;

1.3. use Theorem~\ref{ET} to evaluate $E_{2k}$ for all $q$'s from $T$,
yielding $(\gamma_r)$ for $r = 0, \ldots, \ell+1$\;

1.4 \For{$r:=2$ \KwTo $\ell+1$}{
  \If{$\mathcal{L}_r$ is not solved}{

1.4.1 Instantiate $\mathcal{S}_r$ with $\sum \gamma_i^r$, $A_\rho$ and
  $B_\rho$; add it to $\mathcal{L}_r$\;

1.4.2 \If{$\mathcal{L}_r$ has as many equations as unknowns}{
   solve $\mathcal{L}_r$ and store the values; declare $\mathcal{L}_r$ solved\;
  }

  }
}
}
2. Round the coefficients and use Newton's formulas.
}
\end{algorithm}
The system $\mathcal{S}_r$ has size $N_{2, 3}(r)^2 = O(r^2)$ and we
need $O(r^\omega)$ operations to solve it, for a total of
$O(\ell^{\omega+1})$. Like in the series case, we
can use the $(E_4, E_6, \Delta)$ basis since we anticipate integer
coefficients, as indicated by the series computations.

\smallskip
\noindent
{\bf Example.} Take again $\ell=5$, for which the $\mathcal{S}_r$ were
already given. Let us concentrate on the case of
$\sigma_6(q) = u_{6, 0} E_4^3 + u_{6, 1} \Delta$; we start with
$\mathcal{L}_6 = \emptyset$. Using $\rho = 1.1$ leads to
$$\mathcal{L}_6 = \{
1.912407642 u_{6,0}+0.0009726854527956 u_{6,1} $$
$$= 1393450.57337539139\}$$
and the following iteration with $\rho=1.2$ adds
$$\{
1.435895343 u_{6,0}+0.0005247501300701 u_{6,1} $$
$$= 1156054.63606077432\}$$
and the solution of $\mathcal{L}_6$ (rounded to integers) is
$$u_{6, 1} = -534159360, u_{6, 0} = 1000320.$$

%%%%%%%%%% SS
\subsection{Isogeny volcanoes}

The method in \cite{ChLa05} shares many common points
with the method to be described next but
with a worse complexity. It uses supersingular curves whose
complete explicit $\ell$-torsion is required.

The work of \cite{BrLaSu12} is a building block in \cite{Sutherland13}
where direct evaluation of $\Phi_\ell(X, j(E)) \bmod q$ is made
possible using an explicit version of the Chinese remainder theorem
modulo small primes. Our version leads an easy adaptation to this
problem for the CRT polynomials.

%%%%%%%%%%%%%%% SSS
\subsubsection{Quick presentation}

In a nutshell, the algorithm in~\cite{BrLaSu12} performs computations
modulo special primes $p$ satisfying arithmetical conditions: $p \equiv
1\bmod \ell$ and $4 p = t^2 - \ell^2 v^2 D$ in integers $t$ and $v$,
$v$ not a multiple of $\ell$; $D < 0$ is the discriminant of some
auxiliary quadratic field. With these conditions, the so-called class
polynomial $H_D(X)$ splits completely modulo $p$ and its roots are
$j$-invariants of elliptic curves with complex multiplication that are
needed in the algorithms. Basically, the algorithm interpolates data
using the roots of $H_D(X)$. The central point is to compute isogenies
between the two levels of the isogeny volcano associated to $D$ modulo
$p$. We refer the reader to the original article for more properties
related to elliptic curves. For our purpose, we just need to know that
we have isogeny data available and that they can help us computing the
polynomial $U_\ell(X, A, B) \bmod p$ from these data. We refer to the
article for the complexity under GRH, namely $O(\ell^3 (\log \ell)^3
\log\log\ell)$ using $O(\ell^2 \log (\ell p))$ space for suitably
chosen $p$.

%%%%%%%%%%%%%%% SSS
\subsubsection{The algorithm}

In our case, we consider power sums again.
To lighten the exposition, we consider $U_\ell$ only
since the algorithm is the same for $V_\ell$ and $W_\ell$. 

We adapt here a slight modification of the simplified version Algorithm
2.1 of \cite{BrLaSu12} to our needs. All we describe is also valid in
the full version in \cite{BrLaSu12}.
We denote by
$\powersums$ the powersums of $U_\ell$ viewed as a polynomial in $X$
with coefficients in $\Z[Y, Z]$ for $1 \leq r \leq \ell+1$. We can
write $$\powersums_r(Y, Z) = \sum_{2 i_2 + 3 i_3 = r}
c_{r, i_2, i_3} Y^{i_2} Z^{i_3}.$$
These sums will be reconstructed from values
$\powersums_{r, i}(A_i, B_i)$ associated to curves $\EE_i: Y^2 = X^3 +
A_i X+B_i$.

\LinesNotNumbered
\begin{algorithm}[hbt]
\caption{Computing $U_\ell(X, Y, Z) \bmod p$}
\SetKwProg{Fn}{Function}{}{}
\Fn{\ComputeUMod($\ell$, $D$, $H_{D}(X)$, $p$):}{
\Input{$\ell$ an odd prime, $D$ the discriminant of an imaginary
quadratic order $\mathcal{O}$ of discriminant $D$ with class number
$h(D) \geq \ell+2$; $H_{D}$ is the class
polynomials associated to order $\mathcal{O}$; $p$ prime with
$p\equiv 1\bmod\ell$ and $4 p = t^2 - \ell^2 v^2 D$, $v\not\equiv 0
\bmod \ell$}
\Output{$U_\ell(X, Y, Z) \bmod p$}

1. Build the list $\mathcal{J}_D$ containing the roots of $H_{D}(X)$ 
modulo $p$\;

2. For each $j_i \in \mathcal{J}_D$, find a curve $\EE_i$ having
invariant $j_i$ and cardinality $m = p+1-t$; call $\mathcal{C}_D$ this
set of curves\;

3. For each $\EE_i \in \mathcal{C}_D$ find all its neighbors
$\mathcal{N}(\EE_i)$ in the volcano: $2$ horizontal isogenies and
$\ell-1$ on the floor. This is a collection of $\ell+1$ triplets
$(\sigma_1, \As, \Bs)$ obtained via V\'elu's formulas\;

4. For each $\EE_i: Y^2 = X^3 + A_i X + B_i$, compute the power sums
$\powersums_{r, i}(A_i, B_i)$
corresponding to the values of $\sigma_1$ in $(\sigma_1, \As, \Bs) \in
\mathcal{N}(\EE_i)$ using V\'elu's formulas\;

5. Reconstruct the powers sums $\powersums_r(Y, Z) \bmod
p$ using $A$, $B$ and the $\powersums_{r, i}(A_i, B_i)$, solving a
$N_{2, 3}(r) \times N_{2, 3}(r)$ linear system\;

6. \Return{$\powersums$}.
}
\end{algorithm}
Step 3 is done as follows. Select a random point $P$ of order $\ell$
on $E_i/\GFq{p} = [A_i, B_i]$. Compute the rational isogeny $\EE_i
\rightarrow \EE_i' = \EE_i / \langle P\rangle$ using V\'elu's
formulas. If $j(\EE_i')$
is a root of $H_{D}$, then $\EE_i'$ is on the crater and is one of
the two neighbours. If it does not belong to the crater, it belongs to
the floor. Identification details in the general case with multiple
volcanoes is treated in \cite{BrLaSu12}.

For $V_\ell$ (resp. $W_\ell$), replace $\sigma_1$ by $\As$
(resp. $\Bs$) in Step 5 as far as reconstruction is concerned.

%%%%%%%%%%%%%%% SSS
\subsubsection{A numerical example}

Let us give one value for $\ell = 5$. We select $D = -71$ for which
$h(-71) = 7 \geq 5+2$. Consider $p = 1811$. The roots of $H_{-71}(X)$
modulo $p$ are:
$$\mathcal{J}_D = \{313, 1073, 1288, 1312, 1402, 1767, 1808\}.$$
Associated are curves and neighbours for each $j$ value. These can be
found in Table \ref{fig-ell5}. The power sums $\powersums_r$
corresponding to the values are:
$$\begin{array}{r|rrrrrrr}\hline
\EE_i \backslash r & 1 & 2 & 3 & 4 & 5 & 6 \\ \hline
{[1582, 902]} & 0& 105& 1680& 1379& 756& 772 \\
{[1662, 405]} & 0& 527& 1188& 90& 748& 888 \\
{[1451, 1331]} & 0& 1723& 403& 350& 293& 583 \\
{[1013, 747]} & 0& 1133& 18& 1594& 1738& 105 \\
{[224, 753]} & 0& 95& 760& 1790& 1603& 27 \\
{[1128, 1504]} & 0& 155& 669& 1424& 1130& 522 \\
{[91, 725]} & 0& 1793& 1523& 589& 1233& 134 \\
 \hline
\end{array}$$
For instance, $\sigma_6 = c_{030} Y^3 + c_{012} Z^2$, we need
to solve
$$\left\{
 \begin{array}{ccl}
772 &=& c_{030} 1582^3+ c_{012} 902^2 \bmod 1811, \\
888 &=& c_{030} 1662^3+ c_{012} 405^2 \bmod 1811, \\
583 &=& c_{030} 1451^3+ c_{012} 1331^2 \bmod 1811, \\
\cdots &\cdots &\cdots \\
 \end{array}
\right.$$
that is $1565 Y^3 + 1218 Z^2$. The coefficients are:
\begin{eqnarray*}
\sigma_2 & = & 1771 Y \\
\sigma_3& = & 1331 Z \\
\sigma_4& = & 1120 Y^2 \\
\sigma_5& = & 341 Y Z \\
\sigma_6& = & 1565 Y^3 + 1218 Z^2
\end{eqnarray*}

%%%%%% S
\section{Computing $\As$ and $\Bs$ as rational fractions}

From \cite[Theorem 3.9]{NoYaYo20}, there exist polynomials $N_{\ell,
A}$ and $N_{\ell, B}$ of degree less than $\ell+1$ such that
$$\As = \frac{N_{\ell, A}(X, A, B)}{U_\ell'(X)}, \; \Bs =
\frac{N_{\ell, B}(X, A, B)}{U_\ell'(X)}$$
(Only here: $U_\ell'(X) = \frac{\partial U_\ell}{\partial X}$.)
Moreover, $N_{\ell, A}$ (resp. $N_{\ell, B}$) are polynomials with integer
coefficients and of generalized weight $2 \ell+4$
(resp. $2\ell+6$). The authors of the reference use Groebner basis
computations to find the two numerators.

Alternatively, given $U_\ell$, we can start from (the same is true for
$N_{\ell, B}$):
$U_\ell'(X, A, B) \As = N_{\ell, A}(X, A, B)$
and we plug the series to get
$$U_\ell'(\sigma_1(q), A(q), B(q)) (-\ell^4 E_4(q^\ell)) = N_{\ell,
A}(\sigma_1(q), A(q), B(q)).$$
We find the coefficients by solving
a linear system (over $\Q$ or using small primes as already
described). We can precompute the powers of the series for $\sigma_1$,
$A$ and $B$ and remark that $U_\ell'$ and $N_{\ell, A}$ share a lot of
them. Also, the series $E_4(q^\ell)$ is rather sparse, so that the
product with this quantity is fast.

We can also use floating point numbers as above
and recognize integers in the values of the coefficients of $N_{\ell,
A}$.

There is an advantage to compute
$N_{\ell, A}$ and $N_{\ell, B}$ at the same time, sharing as many
powers as possible, all the more in our use of Pippenger's algorithm.
Numerical examples are given in the appendix.

%%%%% S
\section{Implementation and numerical results}

A lot of trials were done using {\sc Maple} programs, some of which
were then rewritten in {\sc Magma} (version 2.26-10), for speed. See
the author's web page. Computing the polynomials for $\ell \leq 100$
takes a few minutes on a classical laptop. Checking them is done using
SEA, as mentioned in~\cite{NoYaYo20}.

We give some examples of the {\em relative height} $\tilde{H}$ for
some of our polynomials. Here $\tilde{H}(P) = H(P)/((\ell+1)\log
\ell)$. Note that these quantities seem to stabilize when $\ell$
increases.

$$\begin{array}{|r|r|r|r|r|r|}\hline
\multicolumn{6}{|c|}{} \\
\ell & \tilde{H}(\Phi_\ell^t) & \tilde{H}(\Phi_\ell^c) &
\tilde{H}(\Phi_\ell^*) & \tilde{H}(U_\ell) & \tilde{H}(U_\ell^*) \\ \hline
5 & 11.243 & 0.762 & -- & 0.526 & -- \\ %% 5
7 & 9.787 & 0.582 & -- & 0.640 & -- \\ %% 7
11 & 10.130 & 1.842 & 1.120 & 0.670 & 0.240 \\ % 11
13 & 9.565 & 0.367 & 0.941 & 0.688 & -- \\ %% 1
17 & 9.581 & 0.958 & 0.714 & 0.690 & -- \\ %% 5
19 & 9.365 & 0.648 & 0.630 & 0.695 & -- \\ %% 7
23 & 9.438 & 1.995 & 0.419 & 0.698 & 0.441 \\ %% 11
\hline
101 & -- & 1.111 & 0.159 & 0.778 & -- \\ %% 5
103 & -- & 0.740 & 0.249 & 0.779 & -- \\ %% 7
107 & -- & 2.218 & 0.228 & 0.781 & 0.493 \\ %% 11
109 & -- & 0.379 & 0.213 & 0.782 & -- \\ %% 1
\hline
\end{array}$$
Data are computed using the polynomials available in {\sc
Magma}: $\Phi_\ell^c$ is called {\em canonical polynomial} and
$\Phi_\ell^*$ is called {\em Atkin polynomial}. In the case of
$\Phi_\ell^c$, the height depends on $\ell \bmod 12$.
See the appendix for more statistics on the sizes of these polynomials.

When $\ell \equiv 11 \bmod 12$, Atkin~\cite{Atkin92b} suggests to
replace $\sigma$ with $f(q) = (\eta(q)\eta(q^\ell))^2$. For instance:
$$U_{11}^*(X)=
X^{12}-990\,\Delta\,X^{6}+440\,\Delta\,E_4\,X^{4}-165\,\Delta\,E_6\,X^{3}$$
$$+22\,\Delta\,E_4^{2}X^{2}-\Delta\,E_4\,E_6\,X-11\,{\Delta}^{2},$$
whose height is smaller than that of the other alternatives.
Using this type of equation for computing isogenies requires more work,
see~\cite{Morain23b}.

%%%%% S
\section{Conclusions}

We have given a lot of methods for computing the Charlap-Coley-Robbins
polynomials, including representations as fractions in
polynomials. In isogeny cryptography they are useful
for relatively small $\ell$'s, if we store $(U_\ell, N_{\ell, A},
N_{\ell, B})$. If one wants to compute an isogeny, it is enough to
compute a root of $U_\ell$ followed by instantiations of three
polynomials.

There is an alternative to this, suggested by Atkin~\cite{Atkin92b},
using partial derivatives of $U_\ell$. This is described
in~\cite{Morain23b}.

\bigskip
\noindent
{\bf Acknowledgments.} The author wishes to thank A.~Bostan and
F.~Chyzak for helpful discussions around some aspects of this work; special
thanks to the former for his impressive list of references for the
fast evaluation of hypergeometric functions. Thanks also to L.~De Feo for
his updates on cryptographic applications of isogenies.

\def\noopsort#1{}\ifx\bibfrench\undefined\def\biling#1#2{#1}\else\def\biling#1#2{#2}\fi\def\Inpreparation{\biling{In
  preparation}{en
  pr{\'e}paration}}\def\Preprint{\biling{Preprint}{pr{\'e}version}}\def\Draft{\biling{Draft}{Manuscrit}}\def\Toappear{\biling{To
  appear}{\`A para\^\i tre}}\def\Inpress{\biling{In press}{Sous
  presse}}\def\Seealso{\biling{See also}{Voir
  {\'e}galement}}\def\Editor{\biling{Ed.}{R{\'e}d.}}

\appendix

%%%%% S
\section{Some values of $U_\ell$, etc.}
\label{sct:UVWNANB}

Note there is a sign flip compared to~\cite{NoYaYo20} due to a normalization
different from the reference, but used in other articles.

For $\ell = 3$, we compute
$$V_3(X, A, B) =
X^4-84 A X^3+246 A^2 X^2+(63756 A^3+432000 B^2) X$$
$$+576081 A^4+3888000 A B^2,$$
$$W_3(X, A, B) =
X^4+732 B X^3+(25088 A^3+171534 B^2) X^2$$
$$+(1630720 A^3 B+11009548 B^3) X-
139150592 A^3 B^2-437245479 B^4$$
$$-\frac{297493504}{27} A^6;$$
$$N_{3, A}(X, A, B) = 84 X^3 A-360 X^2 B-76 X A^2+36 A B,$$
$$N_{3, B}(X, A, B) = 732 X^3 B+\frac{1456}{3} X^2 A^2-724 X A
B-\frac{112}{3} A^3+108 B^2.$$
$$N_{5, A}
= 630 A X^5-9360 B X^4-8240 A^2 X^3+24480 B A X^2$$
$$+
(1120 A^3- 28800   B^2) X- 3200 B A^2,
$$
$$N_{5, B} =
15630\,X^{5}B+34720\,X^{4}A^{2}-208240\,X^{3}AB$$
$$+ \left( -76160
\,A^{3}+110400\,B^{2} \right) X^{2}+138720\,XA^{2}B-83200\,A B^{2}.
$$

For $P$ a polynomial with rational coefficients, we approximate its
size by the sum $S$ of the number of bits of the absolute value of
each coefficient.
For instance, we this criterion, we find that $S(U_5) = 36$. The total
for $(U_5, N_{5, A}, N_{5, B})$ is $36+91+117=244$ compared to
$S(V_5)=349$, $S(W_5)=602$; $S(\Phi_5^t) = 2838$, $S(\Phi_5^c) =
55$. We note $S_{tot}(\ell) = S(U_\ell) + S(N_{\ell, A}) + S(N_{\ell,
B})$.

$$\begin{array}{|r|r|r|r|r|r|}\hline
\multicolumn{6}{|c|}{} \\
\ell & \tilde{H}(V_\ell) & \tilde{H}(W_\ell) & \tilde{H}(N_{\ell, A})
& \tilde{H}(N_{\ell, B}) & S_{tot}(\ell)\\ \hline
5 & 3.266 & 4.336 & 1.063 & 1.268 & 244 \\
7 & 3.050 & 4.207 & 0.973 & 1.167 & 551 \\
11 & 2.939 & 3.979 & 0.896 & 1.016 & 1816 \\
13 & 2.856 & 3.969 & 0.864 & 0.983 & 2684\\
17 & 2.770 & 3.883 & 0.831 & 0.919 & 5762\\
19 & 2.754 & 3.831 & 0.820 & 0.901 & 7866\\
23 & 2.723 & 3.764 & 0.799 & 0.869 & 13632 \\
\hline
101 & 2.471 & 3.527 & 0.801 & 0.818 & 1187388 \\
103 & 2.469 & 3.517 & 0.801 & 0.818 & 1262059 \\
107 & 2.466 & 3.516 & 0.803 & 0.819 & 1422237 \\
109 & 2.467 & 3.515 & 0.803 & 0.819 & 1504765 \\
\hline
\end{array}$$

%%%%% S
\section{Numerical data for the isogeny volcano algorithm}
\label{sct:num}

\begin{table}
$$\begin{array}{|r|c|c|r|}\hline
i & \EE_i = [A_i, B_i]        & \Es_i & \sigma(\Es_i) \\ \hline
1 & [1582, 902] & [594, 422] & 226 \\
&& [1543, 911] & 1542 \\
&& [937, 1244] & 1283 \\
&& [1333, 561] & 1691 \\
&& [879, 342] & 1212 \\
&& [757, 1578] & 1290 \\
\hline
2 & [1662, 405] & [1770, 433] & 529 \\
&& [1439, 1411] & 1536 \\
&& [259, 355] & 1810 \\
&& [382, 1793] & 1733 \\
&& [1472, 543] & 433 \\
&& [413, 1603] & 1203 \\
\hline 
3 & [1451, 1331] & [1096, 1433] & 743 \\
&& [1371, 1367] & 98 \\
&& [1105, 1195] & 207 \\
&& [1657, 1699] & 787 \\
&& [811, 812] & 1769 \\
&& [779, 1311] & 18 \\
\hline
4 & [1013, 747] & [1691, 473] & 1705 \\
&& [509, 342] & 1245 \\
&& [1642, 417] & 1406 \\
&& [127, 765] & 1519 \\
&& [905, 1464] & 145 \\
&& [1277, 254] & 1224 \\
\hline
5 & [224, 753] & [1485, 892] & 1566 \\
&& [823, 1106] & 908 \\
&& [397, 1451] & 1729 \\
&& [131, 673] & 450 \\
&& [654, 1798] & 1353 \\
&& [1805, 1025] & 1238 \\
\hline
6 & [1128, 1504] & [1275, 1672] & 1176 \\
&& [1409, 761] & 1362 \\
&& [907, 1757] & 309 \\
&& [824, 1267] & 781 \\
&& [578, 1320] & 1208 \\
&& [1168, 1207] & 597 \\
\hline
7 & [91, 725] & [1184, 542] & 1284 \\
&& [1753, 297] & 859 \\
&& [1440, 1524] & 1268 \\
&& [421, 410] & 517 \\
&& [1626, 1013] & 245 \\
&& [198, 159] & 1260 \\
\hline
\end{array}$$
\caption{Values for $\ell = 5$ and $p = 1811$. \label{fig-ell5}}
\end{table}

\end{document}